\documentclass[smallextended,envcountsect]{svjour3}
\smartqed
\usepackage{graphicx}
\journalname{JOTA}
\usepackage{amsfonts}
\usepackage{amssymb}
\usepackage{graphicx}
\usepackage{cancel}
\usepackage{tikz}

\def\hat{\widehat}

\def\emp{\emptyset}
\def\conv{{\rm conv}\,}

\def\dom{{\rm dom}\,}
\def\span{{\rm span}\,}
\def\epi{{\rm epi\,}}

\def\N{{\cal N}}

\def\O{{\cal O}}

\def\sub{\partial}
\def\B{\mathbb B}

\def\ox{\overline{x}}
\def\oy{\overline{y}}
\def\oz{\overline{z}}

\def\disp{\displaystyle}

\def\tto{\;{\lower 1pt\hbox{$\rightarrow$}}\kern-10pt
\hbox{\raise 2pt\hbox{$\rightarrow$}}\;}

\def\Tilde{\widetilde}
\def\Bar{\overline}
\def\ra{\rangle}
\def\la{\langle}

\def\B{I\!\!B}

\def\R{\mathbb{R}}
\def\N{I\!\!N}
\def\ox{\bar{x}}
\def\oy{\bar{y}}
\def\oz{\bar{z}}
\def\ov{\bar{v}}
\def\ow{\bar{w}}
\def\ou{\bar{u}}

\def\op{\bar{p}}

\def\co{{\rm co}}
\def\cone{{\rm cone}}
\def\ri{\mbox{\rm ri}\,}
\def\int{\mbox{\rm int}\,}
\def\gph{\mbox{\rm gph}\,}
\def\epi{\mbox{\rm epi}\,}
\def\dim{\mbox{\rm dim}\,}
\def\dom{\mbox{\rm dom}\,}

\def\ker{\mbox{\rm ker}\,}
\def\conv{\mbox{\rm conv}\,}

\def\lip{\mbox{\rm lip}\,}

\def\aff{\mbox{\rm aff}\,}

\def\b{\hfill\Box}

\def\O{\Omega}
\def\ph{\varphi}
\def\emp{\emptyset}
\def\st{\stackrel}
\def\oR{\Bar{\R}}
\def\lm{\lambda}
\def\gg{\gamma}
\def\dd{\delta}
\def\al{\alpha}

\def \N{I\!\!N}
\def\th{\theta}
\def\vt{\vartheta}
\def\beq{\begin{equation}}
\def\eeq{\end{equation}}

\begin{document}
\title{Stability Analysis for Composite Optimization Problems and Parametric Variational Systems}
\author{B. S. Mordukhovich \and M. E. Sarabi}

\institute{Boris S. Mordukhovich, Corresponding author  \at
Wayne State University,  Detroit, Michigan, USA
\at and  Peoples' Friendship University of Russia, Moscow 117198, Russia\\
\email{boris@math.wayne.edu},
\and
M. Ebrahim Sarabi \at
Miami University, Oxford, Ohio, USA\\
\email{sarabim@miamioh.edu}}
\date{Received: date / Accepted: date}

\maketitle

\begin{abstract}
This paper aims  to provide various applications for second-order variational analysis of  extended-real-valued piecewise liner functions recently obtained in \cite{ms15}. We mainly focus here on establishing relationships between full stability of local minimizers in composite optimization and Robinson's strong regularity of associated (linearized and nonlinearized) KKT systems. Finally, we address Lipschitzian stability of parametric variational systems with convex piecewise linear potentials.
\end{abstract}
\keywords{variational analysis and optimization \and piecewise linear functions \and second-order subdifferentials \and nondegeneracy \and full stability of local minimizers\and strong regularity of KKT systems \and Lipschitzian stability of parametric variational systems}
\subclass{49J52 \and 90C30 \and 90C31}

\section{Introduction}

This paper is dedicated to Professor Boris Teodorovich Polyak who is one of the founders of modern optimization theory and is a very active contributor to current hot areas of optimization, control, and their applications. In particular, we refer the reader to Polyak's classical monograph \cite{pol87} and more recent highly cited publication \cite{pol98} somewhat related to the topics of our paper.

In this paper we continue the path initiated in \cite{ms16} about applications of the second-order subdifferential theory for {\em convex piecewise linear} (CPWL) extended-real-valued functions \cite{rw}, which has been recently developed in \cite{ms15}. Employing here explicit calculations of the {\em second-order subdifferentials} (or generalized Hessians) of such functions in the sense of \cite{m92} together with the second-order characterization of fully stable local optimal solutions for composite optimization problems via the {\em composite SSOSC} (strong second-order sufficient condition) under a certain {\em partial nondegeneracy} allows us to establish comprehensive relationships between several notions of stability in composite optimization and parametric variational systems.

The first topic of our applications concerns {\em Karush-Kuhn-Tucker} (KKT) {\em systems} associated with composite optimization problems described by fully amenable compositions involving CPWL functions. We show, by using a reduction approach combined with second-order calculus and subdifferential computation, that the aforementioned SSOSC effectively characterizes Robinson's {\em strong regularity} \cite{rob} of the corresponding KKT systems arising from composite optimization and also an appropriate Lipschitzian counterpart of Kojima's {\em strong stability} \cite{ko} for perturbed composite optimization problems with respect to ${\cal C}^2$-smooth parameterizations.

Finally, the developed second-order theory \cite{ms15,ms16} helps us to efficiently study {\em robust Lipschitzian stability} of set-valued solution maps to {\em parametric variational systems} (PVS), which are described by subdifferential mappings generated by CPWL functions as well as fully amenable compositions. Employing the {\em coderivative criterion} \cite{m93} for the Lipschitz-like property of multifunctions together with the second-order calculus and precise calculations of the second-order subdifferentials allows us to derive here, depending on the assumptions made, either complete characterizations or effective sufficient conditions for this property explicitly formulated via the initial data.

The rest of the paper is organized as follows. Section 2 recalls some definitions and facts from generalized differential theory of variational analysis needed for the formulations and proofs of the subsequent results. Moreover, for the reader's convenience, we overview basic definitions and results from \cite{ms15,ms16} to make the paper fully self-contained.

In Section~3 we begin our stability analysis for composite optimization models involving CPWL functions by concentrating on {\em strong regularity} of the associated KKT systems and (Lipschitzian) {\em strong stability} of the related stationary points under perturbations. We prove that these differently defined notions occur to be {\em equivalent}  to full stability of the corresponding local minimizers under partial nondegeneracy, being therefore completely characterized by the aforementioned SSOSC.

Section~4 concerns not optimization problems but solution maps to parametric variational systems in the form of perturbed {\em generalized equations} with set-valued subdifferential mappings generated by CPWL functions and fully amenable compositions. We present here comprehensive results for {\em Lipschitzian stability} (in the sense of the validity of the Lipschitz-like/Aubin property) of solution maps expressed entirely of the initial PVS data. In the concluding Section~5  we formulate some unsolved problems of our future research.

The notation used are standard in variational analysis; see \cite{rw,m06}.

\section{Basic Definitions and Preliminaries}\label{prel}

Let us first recall some basic constructions of generalized differentiation exploited in what follows. For $\O\subset\R^n$ with $\ox\in\O$, the {\em normal cone} to $\O$ at $\ox$ (known also as the limiting, basic or Mordukhovich normal cone) is
\begin{eqnarray}\label{2.2}
\begin{array}{ll}
N(\ox;\O):=\Big\{v\in\R^n:&\exists\,x_k\to\ox,\;v_k\to x\;\mbox{ as }\;k\to\infty\;\mbox{ such that }\\
&x_k\in\O\;\mbox{ and }\;\disp\limsup_{x\st{\O}\to x_k}\frac{\la v_k,x-x_k\ra}{\|x-x_k\|}\le 0\Big\},
\end{array}
\end{eqnarray}
where $k\in\N:=\{1,2,\ldots\}$. It is well known that, despite the intrinsic nonconvexity of (\ref{2.2}) for nonconvex sets, the normal cone (\ref{2.2}) possess comprehensive calculus rules; see \cite{rw,m06}. For a set-valued mapping $F\colon\R^n\rightrightarrows\R^m$, its domain and graph are given by
$$
\dom F:=\Big\{x\in\R^n \colon\;F(x)\ne\emp\Big\},\quad\gph F:=\Big\{(x,y)\in\R^n\times\R^m \colon\;x\in F(x)\Big\}.
$$

Define now the {\em coderivative} of $F$ at $(\ox,\oy)\in\gph F$  by
\begin{equation}\label{2.8}
D^*F(\ox,\oy)(u):=\Big\{v\in\R^n \colon\;(v,-u)\in N((\ox,\oy);\gph F)\Big\},
\end{equation}
which is an ``adjoint derivative" of set-valued mappings and reduces to the adjoint/transposed Jacobian operator $D^*f(\ox)(u)=\{\nabla f(\ox)^*u\}$, $u\in\R^m$, if $F=f\colon\R^n\to\R^m$ with $\oy=f(\ox)$ is a single-valued smooth mapping. Recall also that $F$ admits a {\em single-valued graphical localization} around $(\ox,\oy)\in\gph F$ if there exist some neighborhoods $U$ of $\ox$ and $V$ of $\oy$ together with a single-valued mapping $f\colon U\to V$ such that $\gph F\cap(U\times V)=\gph f$.

For an extended-real-valued function $\ph\colon\R^n\to\oR:=]-\infty,\infty ]$, the (first-order) {\em subdifferential} of $\ph$ at $\ox\in\dom\ph:=\{x\in\R^n\colon\;\ph(x)<\infty\}$ is defined via the normal cone (\ref{2.2}) to the epigraph $\epi\ph:=\{(x,\al)\in\R^{n+1}\colon\;\al\ge\ph(x)\}$ by
\begin{equation}\label{2.6}
\partial\ph(\ox):=\Big\{v\in\R^n \colon\;(v,-1)\in N((\ox,\ph(\ox));{{\rm\small epi}\,\ph})\Big\}.
\end{equation}
The {\em second-order subdifferential} of $\ph$ at $\ox\in\dom\ph$ relative to $\ov\in\partial\ph(\ox)$ is defined via the following dual ``derivative-of-derivative" approach of \cite{m92} by
\begin{eqnarray}\label{2nd}
\partial^2\ph(\ox,\ov)(u)\colon=(D^*\partial\ph)(\ox,\ov)(u),\quad u\in\R^n.
\end{eqnarray}

It is worth noticing that the second-order construction (\ref{2nd}) corresponds to the Hessian mapping $\partial^2\ph(\ox,\nabla\ph(\ox))(u)=\{\nabla^2\ph(\ox)u\}$ if $\ph$ is ${\cal C}^2$-smooth around $\ox$. The second-order construction (\ref{2nd}) has been well understood for important classes of functions, mostly appeared in the study of stability of constrained optimization problems; see, e.g., \cite{ms15,dr96,hmn} for more details. Below we briefly recall the explicit calculations of (\ref{2nd}), which have been done recently for the class of convex piecewise linear functions $\th\colon\R^m\to\Bar\R$; we use the notation $\th\in CPWL$ to indicate that $\th$ belongs to this class of functions. Employing \cite[Theorem~2.49]{rw}, we can equivalently state that $\th\in CPWL$ if there are $\alpha_i\in\R$, and  $a_i\in\R^m$, $l\in\N$ with $i\in T_1\colon=\{1,\ldots,l\}$ such that $\th$ is represented via the indicator function of the domain of $\th$ by
\begin{equation}\label{theta}
\th(z)=\max\Big\{\la a_1,z\ra-\alpha_1,\ldots,\la a_l,z\ra-\alpha_l\Big\}+\dd(z;\dom\th),\quad z\in\R^m,
\end{equation}
where the domain set $\dom\th$ is a convex polyhedron given by
\begin{equation}\label{dom}
\dom\th=\Big\{z\in\R^m\colon\;\la d_i,z\ra\le\beta_i\;\mbox{ for all }\;i\in T_2:=\{1,\ldots,p\}\Big\}
\end{equation}
with some elements $d_i\in\R^m$, $\beta_i\in\R$, and $p\in\N$. The domain of $\th$ admits by \cite[Proposition~3.2]{ms15} the union representation $\dom\th=\bigcup^{l}_{i=1}{C_i}$ with $l$ taken from (\ref{theta}) and with the sets $C_i$, $i\in T_1$, defined by
\begin{equation}\label{pwlr1}
C_i:=\Big\{z\in\dom\th\colon\;\la a_j,z\ra-\al_j\le\la a_i,z\ra-\al_i,\;\;\mbox{for all}\;\;j\in T_1\Big\}.
\end{equation}
Defining next the active index subsets for $\oz\in\dom\th $ by
\begin{equation}\label{active2}
K(\oz):=\Big\{i\in T_1\colon\;\oz\in C_i\Big\}\;\mbox{ and }\;I(\oz):=\Big\{i\in T_2\colon\;\la d_i,\oz\ra=\beta_i\Big\},
\end{equation}
we obtain from \cite[Proposition~3.3]{ms15} the following formula for $\partial\th(\oz)$:
\begin{equation}\label{fos}
\partial\th(\oz)=\conv\Big\{a_i:i\in K(\oz)\Big\}+\cone\Big\{d_i\colon\;i\in I(\oz)\Big\}.
\end{equation}
Picking $(\oz,\ov)\in\gph\partial\th$, we conclude from (\ref{fos}) that $\ov=\ov_1+\ov_2$, where
\begin{equation}\label{eq06}
\begin{array}{lll}
\disp\ov_1=\sum_{i\in K(\oz)}\bar\lm_i a_i&\mbox{ with }&\disp\sum_{i\in K(\oz)}\bar\lm_i=1,\;\bar\lm_i\ge 0,\;\mbox{ and }\\
\disp\ov_2=\sum_{i\in I(\oz)}\bar\mu_id_i&\mbox{ with }&\bar\mu_i\ge 0.
\end{array}
\end{equation}
Taking into account representation (\ref{eq06}), we define the two index subsets of positive multipliers by
\begin{equation}\label{eq05}
J_+(\oz,\ov_1):=\Big\{i\in K(\oz)\colon\;\bar\lm_i>0\Big\},\quad J_+(\oz,\ov_2):=\Big\{i\in I(\oz)\colon\;\bar\mu_i>0\Big\}.
\end{equation}
Select arbitrary index subsets $P_1\subset Q_1\subset T_1$ and $P_2\subset Q_2\subset T_2$ and consider the following sets defined entirely via the parameters in (\ref{theta}) and (\ref{dom}):
\begin{eqnarray}\label{eq080}
\begin{array}{lll}
{\cal F}_{\tiny\{P_1,Q_1\},\{P_2,Q_2\}}:&=\span\Big\{a_i-a_j\,\colon\;i,j\in P_1\Big\}\\
&+\cone\Big\{a_i-a_j\,\colon\;(i,j)\in(Q_1\setminus P_1)\times P_1\Big\}\\
&+\span\Big\{d_i\,\colon\;i\in P_2\Big\}+\cone\Big\{d_i\,\colon\;i\in Q_2\setminus P_2\Big\},
\end{array}
\end{eqnarray}
\begin{eqnarray}\label{eq081}
\begin{array}{ll}
{\cal G}_{\tiny\{P_1,Q_1\},\{P_2,Q_2\}}:=\Big\{u\in\R^n\colon&\la a_i-a_j,u\ra=0\;\mbox{if}\;i,j\in P_1,\\
&\la a_i-a_j,u\ra\le 0\;\mbox{if}\;(i,j)\in(Q_1\setminus P_1)\times P_1,\\
&\la d_i,u\ra=0\;\mbox{if}\;i\in P_2,\\
&\la d_i,u\ra\le 0\;\mbox{if}\;i\in Q_2\setminus P_2\;\Big\}.
\end{array}
\end{eqnarray}
Appealing to the above constructions for $\th\in CPWL$ given by (\ref{theta}), we deduce from \cite[Theorem~5.1]{ms15}  the following precise calculation formulas for the second-order subdifferential at any $u\in\R^m$:
\begin{equation}\label{2nd-val}
\begin{array}{ll}
\partial^2\th(\oz,\ov)(u)=\Big\{w\,\colon &(w,-u)\in{\cal F}_{\tiny\{P_1,Q_1\},\{P_2,Q_2\}}\times{\cal G}_{\tiny\{P_1,Q_1\},\{P_2,Q_2\}},\\
&(P_1,Q_1,P_2,Q_2)\in{\cal A}\Big\},
\end{array}
\end{equation}
where the set ${\cal A}$ of index quadruples is defined by
\begin{eqnarray}\label{eq092}
\begin{array}{ll}
{\cal A}:=\Big\{(P_1,Q_1,P_2,Q_2)\,\colon&P_1\subset Q_1\subset K,\;P_2\subset Q_2\subset I,\\
&(P_1,P_2)\in D(\oz,\ov),\;H_{\tiny\{Q_1,Q_2\}}\ne\emp\Big\}
\end{array}
\end{eqnarray}
with $K:=K(\oz)$, $I:=I(\oz)$, $H_{\tiny\{Q_1,Q_2\}}:=\{z\in\dom\th\colon K(z)=Q_1,I(z)=Q_2\}$,
\begin{eqnarray*}
D(\oz,\ov):=\Big\{(P_1,P_2)\subset K\times I\,\colon\;\ov\in\co\{a_i|\;i\in P_1\}+\cone\{d_i\,\colon\;i\in P_2\}\Big\}.
\end{eqnarray*}
Moreover, the second-order subdifferential analysis in \cite[Theorem~5.2]{ms15} gives us the domain formula
\begin{equation}\label{domcod}
\begin{array}{ll}
\dom\sub^2\th(\oz,\ov)=\Big\{u\,\colon& \la a_i-a_j,u\ra=0\;\mbox{ for }\;i,j\in\Gamma(J_1),\\&\la d_t,u\ra=0\;\mbox{ for }\;t\in\Gamma(J_2)\Big\},
\end{array}
\end{equation}
where the index sets $\Gamma(J_1)$ and $\Gamma(J_2)$ are defined by
\begin{equation}\label{feature}
\begin{array}{ll}
\Gamma(J_1):=\Big\{i\in K\colon\la a_i-a_j,u\ra=0\;\mbox{for}\;j\in J_1,\;u\in{\cal G}_{\tiny\{J_1,K\},\{J_2,I\}}\;\Big\},\\
\Gamma(J_2):=\Big\{t\in I\colon\la d_t,u\ra=0\;\mbox{for}\;u\in{\cal G}_{\tiny\{J_1,K\},\{J_2,I\}}\;\Big\}
\end{array}
\end{equation}
with the notation $J_1:=J_+(\oz,\ov_1)$ and $J_2:=J_+(\oz,\ov_2)$ from (\ref{eq06}) and (\ref{eq05}).

Given $\th\in CPWL$, we always assume in what follows that $0\in \aff \sub\th(\oz)$ with $\aff\sub\th(\oz)$ being the affine hull of $\sub\th(\oz)$. As it was explained in \cite[Section~3]{ms16}, there is no loss of generality in assuming the latter condition when we deal with the second-order subdifferential (\ref{2nd}). Indeed, we have $S(\oz)=\aff\partial\theta(\oz)-b_{\oz}$ for some $b_{\oz}\in\aff\partial\theta(\oz)$. Defining then $\bar{\theta}(z):=\theta(z)-\la b_{\bar{z}},z\ra$ shows that $0\in\aff\partial\bar\th(z)$ and $\partial^2\th(\oz,\oy)=\partial^2\bar\th(\oz,\oy-b_{\oz})$ for any $\ov\in\partial\th(\oz)$.
\textcolor{red}{\xcancel{This} Employing  $0\in \aff \sub\th(\oz)$} leads us to
\begin{equation}\label{ss1}
S(\oz)=\aff\sub\th(\oz),
\end{equation}
where $S(\oz)$ stands for a subspace of $\R^m$ parallel to the affine hull $\aff\sub\th(\oz)$. It is proved in \cite[Lemma~3.1]{ms16} that such a CPWL function $\th\colon\R^m\to\oR$ is {\em ${\cal C}^\infty$-reducible} to a function $\vt\colon\R^s\to\oR$ at $\oz$ with $s=\dim S(\oz)\le m$; this means that  there exists a ${\cal C}^\infty$-smooth  mapping $h\colon\R^m\to\R^s$ with the surjective derivative $\nabla h(\oz)$ such that $\th(z)=(\vartheta\circ h)(z)$ for all $z$ around $\oz$. The interested readers can find more about reducibility of functions and sets in \cite{ms16,mnn,bs}.

Following \cite{mnn}, we say that a pair $(\ox,\ow)\in\R^n\times\R^d$ is a {\em partial nondegenerate point} of $\Phi\colon\R^n\times\R^d\to\R^m$ in $x$ relative to the mapping $h\colon\R^m\to\R^s$ taken from the ${\cal C}^\infty$-reducibility of $\th$ if
\begin{equation}\label{fnond}
\nabla_x\Phi(\ox,\ow)\R^n+\ker\nabla h(\oz)=\R^m\;\mbox{ with }\;\oz=\Phi(\ox,\ow).
\end{equation}
The above definition of the nondegenerate points is an extension of the nondegeneracy concept for sets; see \cite[Definition~4.70]{bs} for more details. It is important to point out that for standard problems of nonlinear programming with smooth data this concept reduces to the classical linear independence constraint qualification (LICQ) as observed in \cite[Example~4.77]{bs}. In the subsequent sections of the paper we refer to (\ref{fnond}) as the nondegeneracy condition (ND). It has been recently observed in \cite[Theorem~3.2]{ms16} that condition ND admits the following dual equivalent representation:
\begin{equation}\label{2qcf2}
S(\oz)\cap\ker\nabla_x\Phi(\ox,\ow)^*=\{0\}
\end{equation}
with $S(\oz)$ in (\ref{ss1}). Note \textcolor{red}{\xcancel{we} {that}} the explicit calculation of $S(\oz)$ entirely via the initial parameter of the CPWL function $\theta$ is given in \cite[Theorem~3.1]{ms16}.

Observe to this end that the inner mapping $h$ in the definition of reducibility is not necessarily unique. Moreover,the reducibility  is always used together with the nondegeneracy condition. This fact brings some limitations to the kind of $h$ that can be selected in the definition of reducibility. For instance, if we choose $h=Id$, i.e., $\th$ is reducible to itself, then (\ref{fnond}) tells us that
$$
\nabla_x\Phi(\ox,\ow)\R^n+\ker\nabla h(\oz)=\nabla_x\Phi(\ox,\ow)\R^n=\R^m,
$$
which says that the operator $\nabla_x\Phi(\ox,\ow)$ is surjective. However, it is shown in \cite{ms16} that we can find a mapping $h$ for which the nondegeneracy condition (\ref{fnond}) offers a strictly weaker constraint qualification.

In this paper we often deal with a composition $\th\circ\Phi$ of CPWL outer functions $\th\colon\R^m\to\oR$, and inner mappings $\Phi\colon\R^n\times\R^d\to\R^m$ that are ${\cal C}^2$-smooth around some $(\ox,\ow)$ with $\oz:=\Phi(\ox,\ow)\in\dom\th$ under the following {\em first-order qualification condition}:
\begin{equation}\label{2.9}
\sub^{\infty}\th(\oz)\cap\ker\nabla_x\Phi(\ox,\ow)^*=\{0\}.
\end{equation}
Such compositions are an important subclass of functions called {\em fully amenable}   in $x$ at $\ox$ with compatible parametrization by $w$ at $\ow$; see \cite{rw,lpr} for more details.

\section{Strong Regularity and Strong Stability in Composite Models}\label{str-reg}

In this section we mainly focus on the study of stability analysis of the composite optimization problem given by
\begin{equation}\label{fcp1}
\mbox{minimize }\ph_0(x)+\theta\big(\Phi(x)\big)\;\mbox{s.t.}\;x\in\R^n\;\mbox{with}\;\Phi(x):=\big(\ph_1(x),\ldots,\ph_m(x)\big),
\end{equation}
where $\theta\colon\R^m\to\Bar\R$ is a CPWL extended-real-valued function, and where all $\ph_i\colon\R^n\to\R$, $i=0,\ldots,m$, are ${\cal C}^2$-smooth around the reference optimal solution. As argued in \cite{ms16}, the presented model provides a very convenient framework for developing theoretical aspects of optimization in broad classes of constrained problems including nonlinear programs as well as constrained and unconstrained minimax problems. Note that, besides the aforementioned classes, the composite optimization format under consideration includes the following major subclass of {\em extended nonlinear programs} given by
$$
\mbox{minimize }\;\ph_0(x)+(\th\circ\Phi)(x)\;\mbox{ with }\;\th(x):=\disp\sup_{p\in P}\la p,x\ra,\quad x\in\R^n,
$$
where $P\subset\R^n$ is a convex polyhedra, and where $\th\colon\R^m\to\oR$ is CPWL.

To proceed, consider the two-parametric version of (\ref{fcp1}) constructed by
\begin{equation}\label{fcp2}
{\cal P}(w,v):\quad\mbox{minimize }\;\ph_0(x,w)+\theta(\Phi(x,w))-\la v,x\ra\;\mbox{ subject to }\;x\in\R^n,
\end{equation}
where the perturbed functions $\ph_0(x,w)$ and $\Phi(x,w)=(\ph_1(x,w),\ldots,\ph_m(x,w))$ are ${\cal C}^2$-smooth with respect to both variables. Denote
\begin{equation}\label{phi}
\ph(x,w):=\ph_0(x,w)+\th(\Phi(x,w))\;\mbox{ for }\;(x,w)\in\R^n\times\R^d
\end{equation}
and then fix $\gg>0$ and $(\ox,\ow,\ov)$ with $\Phi(\ox,\ow)\in\dom\th$ and $\ov\in\partial_x\ph(\ox,\ow)$. Define the parameter-depended  optimal value function for (\ref{fcp2}) by
\begin{eqnarray*}
m_\gg(w,v):=\inf_{\|x-\ox\|\le\gg}\Big\{\ph(x,w)-\la v,x\ra\Big\}
\end{eqnarray*}
and the parameterized set of optimal solutions to (\ref{fcp1}) by
\begin{eqnarray}\label{M}
M_\gg(w,v):=\mbox{argmin}_{\|x-\ox\|\le\gg}\Big\{\ph(x,w)-\la v,x\ra\Big\}
\end{eqnarray}
with the convention that argmin:=$\emp$ when the expression under minimization is $\infty$. A point $\ox$ is called \cite[Definition~1.1]{lpr} a {\em fully stable} locally optimal solution to problem ${\cal P}(\ow,\ov)$ in (\ref{fcp2}) if there exist a number $\gg>0$ and neighborhoods $W$ of $\ow$ and $V$ of $\ov$ such that the mapping $(w,v)\mapsto M_\gg(w,v)$ is single-valued and Lipschitz continuous with $M_\gg(\ow,\ov)=\{\ox\}$ and the function $(w,v)\mapsto m_\gg(w,v)$ is likewise Lipschitz continuous on $W\times V$.

The seminal notion of full stability of local optimal solutions was first introduced in \cite{lpr} in the extended-real-valued format of unconstrained optimization. Moreover, it was characterized in \cite[Theorem~2.3]{lpr} using  the second-order subdifferential given by (\ref{2nd}). Recently, second-order characterizations of full stability have been established for constrained optimization problems including NLPs, SOCPs, and optimal control of semilinear PDEs; see \cite{ms15,ms16,mn3,mnr,mn15,mrs,mos} for more details and discussions.
Quite recently \cite[Theorem~4.1]{ms16}, we have obtained a second-order characterization of full stability for the composite optimization problem (\ref{fcp2}) entirely in terms of the initial data of (\ref{fcp2}). In this section, employing the latter characterization, we continue our study of {\em composite optimization} problems of type (\ref{fcp1}) with CPWL outer functions $\th$. Our main goal is to establish relationships between full stability of local minimizers in (\ref{fcp2}) and some other stability/regularity notions for perturbed versions of (\ref{fcp1}) and associated (linearized and nonlinearized) KKT systems. The notions under consideration revolve around Robinson's {\em strong regularity} \cite{rob} and the Lipschitzian version of Kojima's {\em strong stability} \cite{ko}. Involving the nondegeneracy condition ND from (\ref{fnond}) in composite optimization together with employing the reduction approach, discussed below, we show that these notions are actually equivalent in our setting while being also equivalent to full stability of local minimizers under appropriate choices of perturbations.

To proceed, let $\oz:=\Phi(\ox,\ow)\in\dom\th$ with $\th\in CPWL$. Remember from Section 2 that $\th$ is  ${\cal C}^\infty$-reducible to a function $\vt\colon\R^s\to\oR$ at $\oz$ with $s\le m$ by a ${\cal C}^\infty$-smooth  mapping $h\colon\R^m\to\R^s$ with the surjective derivative $\nabla h(\oz)$. Moreover, it is proved in \cite[Lemma~3.1]{ms16} that the mapping $h$ has a linear representation $h(z)=Bz$ with $B$ being a $s\times m$ matrix; see the latter lemma for more details about the matrix $B$. This tells us that
\begin{equation}\label{re1}
\th(z)=(\vartheta\circ B)(z) \quad\quad \mbox{for all}\quad z\;\; \mbox{close to}\;\;\oz.
\end{equation}
Defining now the mapping  $\Psi(x,w):=(B\circ\Phi)(x,w)$, $(x,w)\in\R^n\times \R^d$, allows us to get the reduced problem
\begin{equation}\label{fcp3}
{\cal P}_r(w,v):\quad\mbox{minimize }\;\ph_0(x,w)+\vartheta(\Psi(x,w))-\la v,x\ra\;\mbox{ subject to }\;x\in\R^n.
\end{equation}
The reduced problem (\ref{fcp3})  plays an extremely important role in the proofs of the results established in what follows. The main characteristic feature of this problem, obtained in \cite[Proposition~4.1]{ms16}, is that the derivative $\nabla_x\Psi(\ox,\ow)$ has full rank provided that the nondegeneracy condition (\ref{fnond}) holds. Moreover, it is proved in \cite[Proposition~4.1]{ms16} that the full stability of local optimal solutions to problem ${\cal P}(w,v)$ is {\em equivalent} to that of local optimal solutions to ${\cal P}_r(w,v)$.

Taking into account the qualification condition (\ref{2.9}), the stationary condition $\ov\in\partial_x\ph(\ox,\ow)$ via $\ph$ from (\ref{phi}) can be equivalently written as
\begin{equation}\label{ov}
\ov\in\nabla_x\ph_0(\ox,\ow)+\nabla_x\Phi(\ox,\ow)^*\sub\th(\Phi(\ox,\ow)).
\end{equation}
Thus the KKT systems for ${\cal P}(w,v)$ and ${\cal P}_r(w,v)$ are given, respectively, by
\begin{eqnarray}\label{kkt}
\left\{\begin{array}{ll}
v=\nabla_xL(x,w,\lm),\quad\lm\in\sub\th(\Phi(x,w))\\
\mbox{with }\;L(x,w,\lm):=\ph_0(x,w)+\langle\lm,\Phi(x,w)\rangle,
\end{array}
\right.
\end{eqnarray}\vspace*{-0.2in}
\begin{eqnarray}\label{rkkt}
\left\{\begin{array}{ll}
v=\nabla_x {L}_r(x,w,\mu),\quad\mu\in\sub\vartheta(\Psi(x,w))\\
\mbox{with }\;{L}_r(x,w,\mu):=\ph_0(x,w)+\langle\mu,\Psi(x,w)\rangle.
\end{array}
\right.
\end{eqnarray}
It is worth noticing that representation (\ref{re1}) together with the full rank property of the matrix $B$, coming from the surjectivity of $\nabla h(\oz)$, implies that the Lagrange multipliers $\lm$ of (\ref{kkt}) and $\mu$ of (\ref{rkkt}) are related by $\lm=B^*\mu$.

It is often more convenient in what follows to rewrite the KKT system (\ref{kkt}) as the {\em generalized equation}
\begin{equation}\label{kkt11}
\left[\begin{array}{c}
v\\
0\\
\end{array}
\right]
\in\left[\begin{array}{c}
\nabla_x L(x,w,\lm)\\
-\Phi(x,w)\\
\end{array}
\right]+\left[\begin{array}{c}
0\\
({\sub\th})^{-1}(\lm)\\
\end{array}
\right]
\end{equation}
and denote by $S_{KKT}\colon(w,v)\mapsto(x,\lm)$ the solution map to (\ref{kkt11}).

To provide a second-order characterization of fully stable local optimal solutions to problem ${\cal P}(w,v)$, the following composite strong second-order sufficient condition (SSOSC)
was introduced in \cite[Definition~4.1]{ms16}: we say that the {\em composite SSOSC} holds at  $(\ox,\ow,\ov,\bar\lm)\in\R^n\times\R^d\times\R^n\times\R^m$ with $\ov$ and $\bar\lm$ satisfying {\rm(\ref{ov})} and {\rm(\ref{kkt})}, respectively, if
\begin{equation}\label{sssoc1}
\la u,\nabla^2_{xx}L(\ox,\ow,\bar\lm)u\ra>0\;\mbox{ for all }\;0\ne u\in {\cal S},
\end{equation}
where $L$ is the Lagrangian from {\rm(\ref{kkt})} while the subspace ${\cal S}$ is defined by
\begin{equation}\label{sssoc2}
\begin{array}{ll}
{\cal S}:=\Big\{u\in\R^n\,\colon&\la a_i-a_j,\nabla_x\Phi(\ox,\ow)u\ra=0\;\;\mbox{for}\;\;i,j\in\Gamma(J_1),\\
&\la d_t,\nabla_x\Phi(\ox,\ow)u\ra=0\;\;\mbox{for}\;\;t\in\Gamma(J_2)\Big\}
\end{array}
\end{equation}
via the index sets $\Gamma(J_1)$ and $\Gamma(J_2)$ taken from {\rm(\ref{feature})}. As discussed in \cite{ms16}, the composite SSOSC is an adaptation of Robinson's SSOSC \cite{rob}, introduced for classical NLPs, for the composite optimization problem (\ref{fcp2}). Below, we recall from \cite[Theorem~4.1]{ms16} the second-order characterization of full stable local optimal solutions to ${\cal P}(\ow,\ov)$ via the presented composite SSOSC.

\begin{theorem}{\bf(second-order characterization of full stability in composite optimization).}\label{ssooc} Let $\ox$ be a feasible solution to ${\cal P}(\ow,\ov)$ from {\rm(\ref{fcp2})} for the parameter pair $(\ow,\ov)$ with $\ov$ from {\rm (\ref{ov})}, let $\th\in CPWL$, and let $(\oz,\ov)\in\gph\sub\th$ with $\oz=\Phi(\ox,\ov)$. Under the validity of condition ND from {\rm(\ref{fnond})}, let $\bar\lm$ be a unique solution of the KKT system {\rm (\ref{kkt})}. Then $\ox$ is a fully stable local minimizer of ${\cal P}(\ow,\ov)$ if and only if the composite SSOSC from {\rm(\ref{sssoc1})} is satisfied.
\end{theorem}

Robinson's idea \cite{rob} to define the property of strong regularity for generalized equations involved considering Lipschitzian single-valued localizations of solution maps to an appropriate {\em linearization}. This idea was further developed and applied in many publications; see, e.g., \cite{dr,fp,kk} and the references therein. We keep such a definition of strong regularity in the case of (\ref{kkt11}) and study it in this section. However, it is more convenient for us to start with a similar property for the solution map $S_{KKK}$ of the KKT system (\ref{kkt11}) {\em itself}, without any linearization, and characterize it via the composite SSOSC.

\begin{definition}{\bf(SVLL property of KKT systems).}\label{svll} We say that the KKT system {\rm(\ref{kkt11})} associated with the composite optimization problem {\rm(\ref{fcp2})} has the {\sc single-valued Lipschitzian localization (SVLL)} property at $(\ox,\bar\lm,\ow,\ov)\in\gph S_{KKT}$ if its solution map $S_{KKT}\colon(w,v)\mapsto(x,\lm)$ admits a Lipschitzian single-valued graphical localization around $(\ow,\ov,\ox,\bar\lm)$.
\end{definition}

The next theorem shows the SVLL property of (\ref{kkt11}) is characterized by the simultaneous fulfillment of the composite SSOSC and the nondegeneracy condition ND in composite optimization. It extends the corresponding result of \cite[Theorem~4.10]{bs} and \cite[Theorem~6]{dr96} for NLPs; see also commentaries in \cite{dr,kk,bs} on related developments in this direction.

\begin{theorem}{\bf (characterization of SVLL property via ND and composite SSOSC).}\label{rsr} Let $\ox$ be a feasible solution to problem ${\cal P}(\ow,\ov)$ in {\rm(\ref{fcp2})} with some $\ow\in\R^d$ and $\ov$ from $(\rm{\ref{ov}})$, where $\th\in CPWL$ and $\Phi$ is ${\cal C}^2$-smooth around $(\ox,\ow)$. Consider the following statements:

{\bf (i)} The SVLL property from Definition~{\rm\ref{svll}} holds and we have $\ox\in M_\gamma(\ow,\ov)$ for the argminimum set {\rm(\ref{M})} with some $\gg>0$.

{\bf (ii)} Both the composite SSOSC property {\rm(\ref{sssoc1})} and the nondegeneracy condition ND from {\rm(\ref{fnond})} hold.\\[1ex]
Then we have {\rm(ii)}$\Longrightarrow${\rm(i)}, while the converse application is fulfilled if in addition the first-order qualification condition {\rm(\ref{2.9})} is satisfied.
\end{theorem}
{\it Proof} Suppose first that (ii) holds and deduce from Theorem~\ref{ssooc} that $\ox$ is a fully stable locally optimal solution to ${\cal P}(\ow,\ov)$. It follows from \cite[Proposition~4.1]{ms16} that $\ox$ is also a fully stable locally optimal solution to the reduced problem ${\cal P}_r(\ow,\ov)$. Similarly to (\ref{kkt11}), we can write the KKT system for the reduced problem (\ref{rkkt}) in the generalized equation form
\begin{equation}\label{kkt11r}
\left[\begin{array}{c}
v\\
0\\
\end{array}
\right]
\in\left[\begin{array}{c}
\nabla_x L_r(x,w,\mu)\\
-\Psi(x,w)\\
\end{array}
\right]+\left[\begin{array}{c}
0\\
({\sub\vartheta})^{-1}(\mu)\\
\end{array}
\right]
\end{equation}
and denote by $S_{KKT}^r\colon(w,v)\mapsto(x,\mu)$ its solution map. By (\ref{re1}) we have that the representation $\th=\vartheta\circ B$ holds locally around $\oz$. Moreover, it follows from \cite[Lemma~3.1]{ms16} that $\vartheta\in CPWL$. Remembering that $\lm=B^*\mu$, we split the proof of (ii)$\Longrightarrow$(i) into several steps.\vspace{0.05in}\\
{\bf{Step~1:}} {\em The conditions in {\rm (ii)} imply that the solution map for the generalized equation {\rm(\ref{kkt11r})}, denoted by $S_{KKT}^r\colon(w,v)\mapsto(x,\mu)$, has the SVLL property around $(\ow,\ov,\ox,\bar\mu)$.}\\[1ex]
We start the proof of this fact by recalling that the full stability of $\ox$ in ${\cal P}_r(\ow,\ov)$ ensures by \cite[Theorem~3.4]{mrs} that the set-valued mapping
$$
S_r(w,v):=\Big\{x\in\R^n\,\colon\;v\in\nabla_x\ph_0(x,w)+\nabla_x\Psi(x,w)^*\sub\vartheta(\Psi(x,w))\Big\}
$$
admits a Lipschitzian single-valued graphical localization around $(\ow,\ov,\ox)$. Employing this together with the surjectivity of $\nabla_x\Psi(\ox,\ow)$, which comes from the second part of \cite[Proposition~4.1]{ms16}, yields the mapping $S_{KKT}^r\colon(w,v)\mapsto(x,\mu)$ is single-valued around $(\ow,\ov,\ox,\bar\mu)$. Observe that the Lipschitz continuity of $(w,v)\mapsto x_{wv}=:x$ around $(\ow,\ov)$ is a direct consequence of the full stability of $\ox$ in the reduced problem ${\cal P}_r(\ow,\ov)$. Let the latter property hold in some neighborhoods $W$ of $\ow$ and $V$ of $\ov$. To verify the same property for the mapping $(w,v)\mapsto\mu_{wv}=:\mu$, pick $w_1,w_2\in W$ and $v_1,v_2\in V$ and thus find $\mu_{w_iv_i}\in \sub\vartheta(c_i)$ with $c_i:=\Psi(x_{w_iv_i},w_i)$ for $i=1,2$ satisfying
$$
\left\{\begin{array}{l}
v_2=\nabla_x\ph_0(x_{w_2v_2},w_2)+\nabla_x\Psi(x_{w_2v_2},w_2)^*\mu_{w_2v_2},\\
v_1=\nabla_x\ph_0(x_{w_1v_1},w_1)+\nabla_x\Psi(x_{w_1v_1},w_1)^*\mu_{w_1v_1}.
\end{array}\right.
$$
This allows us to obtain the equality
$$\begin{array}{lll}
\nabla_x\Psi(x_{w_2v_2},w_2)^*(\mu_{w_2v_2}-\mu_{w_1v_1})&=&\Big(\nabla_x\Psi(x_{w_1v_1},w_1)-\nabla_x \Psi(x_{w_2v_2},w_2)\Big)^*\mu_{w_1v_1}\\
&+&\nabla_x\varphi_0(x_{w_1v_1},w_1)-\nabla_x\varphi_0(x_{w_2v_2},w_2)+v_2-v_1,
\end{array}
$$
where $\nabla_x\Psi(x_{w_iv_i},w_i)$ are surjective due to this property for $\nabla_x\Psi(\ox,\ow)$. By \cite[Lemma~1.18]{m06} there is ${\kappa}_{wv}>0$ such that
\begin{eqnarray*}
\begin{array}{lll}
\|\nabla_x\Psi(x_{w_2v_2},w_2)^*(\mu_{w_2v_2}-\mu_{w_1v_1})\|&\ge\kappa_{w_2v_2}\|\mu_{w_2v_2}-\mu_{w_1v_1}\|\\&\ge\kappa\|\mu_{w_2v_2}-\mu_{w_1v_1}\|
\end{array}
\end{eqnarray*}
for $(w,v)\in W\times V$, where $\kappa:=\inf\{\kappa_{wv}|\;(w,v)\in W\times V\}$. Now we claim that $\kappa>0$. Indeed, assuming that $\kappa=0$ gives us $(w_k,v_k)\to(\ow,\ov)$ with $\kappa_{w_k v_k}\to 0$. Appealing to \cite[Lemma~1.18]{m06}, we deduce that
 $$\kappa_{w_k v_k}=\inf\Big\{\|\nabla_x\Psi(x_{w_k v_k},w_k)^*y\|\;:\;\|y\|=1\Big\}.$$ This allows us to find $y_k$ with $\|y_k\|=1$ and
\begin{equation}\label{s017}
\|\nabla_x\Psi(x_{w_k v_k},w_k)^*y_k\|<\kappa_{w_k v_k}+\frac{1}{k}.
\end{equation}
Suppose without loss of generality that $y_k\to\oy$ as $k\to\infty$ with $\|\oy\|=1$. Passing to limit in (\ref{s017}) gives us $\nabla_x\Psi(\ox,\ow)^*\oy=0$. Since the derivative  $\nabla_x\Psi(\ox,\ow)$ is surjective due to  \cite[Proposition~4.1]{ms16}, we arrive at $\oy=0$, which is a contradiction telling us that $\kappa>0$. By the surjectivity of  $\nabla_x\Psi(\ox,\ow)$ there exists a constant $\rho<\infty$ so that $\|\mu_{wv}\|\le\rho$ for all $(w,v)\in W\times V$. Denoting by $\ell>0$ a common Lipschitz constant for the mappings $\nabla_x\ph_0$, $\nabla_x\Psi$, and $(w,v)\mapsto x_{wv}$ on $W\times V$ yields
\begin{eqnarray*}
\begin{array}{lll}
\|\mu_{w_2v_2}-\mu_{w_1v_1}\|&\le&{\kappa^{-1}}\Big(\|\nabla_x\Psi(x_{w_1v_1},w_1)-\nabla_x\Psi(x_{w_2v_2},w_2)\|\cdot\|\mu_{w_1v_1}\|\\
&+&\|\nabla_x\varphi_0(x_{w_1v_1},w_1)-\nabla_x\varphi_0(x_{w_2v_2},w_2)\|+\|v_2-v_1\|\Big)\\&\le&{\kappa^{-1}}
\Big[\rho\ell\Big(\|x_{w_2v_2}-x_{w_1v_1}\|+\|w_2-w_1\|\Big)\\&+&\ell\Big(\|x_{w_2v_2}-x_{w_1v_1}\|+\|w_2-w_1\|\Big)+\|v_2-v_1\|\Big],
\end{array}
\end{eqnarray*}
which justifies the local Lipschitz continuity of the mapping $(w,v)\mapsto\mu_{wv}$.
\vspace{0.05in}\\
{\bf{Step~2:}} {\em The conditions in {\rm (ii)} imply that the SVLL property of {\rm(\ref{kkt11})} is satisfied at the point $(\ox,\bar\lm,\ow,\ov)$}.\\[1ex]
The assumptions of (ii) ensure by Theorem~\ref{ssooc} that $\ox$ is a fully stable local minimizer of ${\cal P}(\ow,\ov)$, and so \cite[Theorem~3.4]{mrs} tells us that the mapping
\begin{equation}\label{eq0107}
S(w,v):=\Big\{x\in\R^n\,\colon\;v\in\nabla_x\ph_0(x,w)+\nabla_x\Phi(x,w)^*\sub\th(\Phi(x,w))\Big\}
\end{equation}
is single-valued and locally Lipschitzian around $(\ow,\ov,\ox)$. Since condition ND from (\ref{fnond}) holds, the Lagrange multiplier in (\ref{kkt11}) is unique, and hence the mapping $S_{KKT}\colon(w,v)\mapsto(x,\lm)$ is single-valued around $(\ow,\ov,\ox,\bar\lm)$. Further, the Lipschitz continuity of $(w,v)\mapsto x_{wv}=:x$  around $(\ow,\ov)$ follows from the full stability of $\ox$. Taking $W$ and $V$ from Step~1, pick $w_i\in W$ and $v_i\in V$, $i=1,2$. Using the relationship $\lm=B^*\mu$, for each $i$ find a unique multiplier $\mu_{w_iv_i}\in\sub\vartheta(c_i)$ with $c_i:=\Psi(x_{w_iv_i},w_i)$ so that
$\lm_{w_iv_i}:=B^*\mu_{w_iv_i}$. This yields
\begin{eqnarray*}
\begin{array}{lll}
\|\lambda_{w_2v_2}-\lambda_{w_1v_1}\|&=&\|B^*\mu_{w_2v_2}-B^*\mu_{w_1v_1}\|\\
&\le&\|B^*\|\cdot\|\mu_{w_2v_2}-\mu_{w_1v_1}\|,
\end{array}
\end{eqnarray*}
which thus justifies the local Lipschitz continuity of the mapping $(w,v)\mapsto\lm_{wv}$ due to Step~1. This completes the proof of implication (ii)$\Longrightarrow$(i).

To verify (i)$\Longrightarrow$(ii), suppose that the SVLL condition satisfies and pick $\eta\in S(\oz)\cap\ker\nabla_x\Phi(\ox,\ow)^*$ with $\oz=\Phi(\ox,\ow)$, where $S(\oz)$ comes from (\ref{ss1}). Since $S(\oz)=\aff\sub\th(\oz)$ due to (\ref{ss1}), we get $\eta\in\aff\sub\th(\oz)$ and deduce from (\ref{2.9}) that $S_{KKT}(\ow,\ov)=\{(\ox,\bar\lm)\}$ for some $\bar\lm\in\R^m$. If $\bar\lm\in\ri\sub\th(\oz)$ with ``ri" standing for the relative interior of a convex set, then $\bar\lm+t\eta\in\sub\th(\oz)$ for any small $t>0$, which tells us that $(\ox,\bar\lm+t\eta)\in S_{KKT}(\ow,\ov)$. Employing now the single-valuedness of the mapping $S_{KKT}$, we get $\eta=0$, and so by (\ref{2qcf2}) condition ND from (\ref{fnond}) holds in this case. Suppose now that $\bar\lm\not\in\ri\sub\th(\oz)$ and, taking into account that $\ri\sub\th(\oz)\ne\emp$, pick $\eta\in\ri\sub\th(\oz)$. It follows from \cite[Proposition~2.40]{rw} that $\bar\lm+t(\eta-\bar\lm)\in\ri\sub\th(\oz)$ for any $t\in(0,1)$. Letting $v_t:=t\nabla_x\Phi(\ox,\ow)^*(\eta-\bar\lm)$ for small $t>0$ gives us $(\ox,\bar\lm+t(\eta-\bar\lm))\in S_{KKT}(\ow,\ov+v_t)$. Remember that $\bar\lm+t(\eta-\bar\lm)\in\ri\sub\th(\oz)$, which allows us to repeat the above arguments and to justify the validity of ND.

To end the proof, it is not hard to see by SVLL that the mapping $S(w,v)$ in (\ref{eq0107}) is single-valued and locally Lipschitzian around $(\ow,\ov,\ox)$. Remembering that $\ox\in M_\gamma(\ow,\ov)$ in (i) and appealing to \cite[Theorem~3.4]{mrs}, with taking into account that the qualification condition imposed therein follows from the justified ND, tell us that $\ox$ is a fully stable local minimizer of ${\cal P}(\ow,\ov)$. Thus SSOSC holds by Theorem~\ref{ssooc}, and we complete the proof of the theorem. $\qed$\vspace*{0.05in}

Next we proceed with the definition and second-order characterization of Robinson's strong regularity for the KKT system (\ref{kkt11}) associated with problem ${\cal P}(\ow,\ov)$ of composite optimization.

\begin{definition}\label{str}{\bf(strong regularity of KKT in composite optimization).} Let $(\ox,\bar\lm)$ be a solution to {\rm(\ref{kkt11})} for  $(w,v)=(\ow,\ov)$ with $\ov=0$. We say that $(\ox,\bar\lm)$  is {\sc strongly regular} for KKT {\rm (\ref{kkt11})} if the solution map to the linearized system at $(\ox,\bar\lm)$ defined by
\begin{eqnarray*}\label{sr3}
\left[\begin{array}{c}
v_1\\
v_2\\
\end{array}
\right]
\in\left[\begin{array}{c}
\nabla_{xx}^2 L(\ox,\ow,\bar\lm)(x-\ox)+\nabla_x\Phi(\ox,\ow)^*(\lm-\bar\lm)\\
-\Phi(\ox,\ow)-\nabla_x\Phi(\ox,\ow)(x-\ox)\\
\end{array}
\right]+\left[\begin{array}{c}
0\\
({\sub\th})^{-1}(\lm)\\
\end{array}
\right]
\end{eqnarray*}
admits a Lipschitzian single-valued graphical localization around $(0,0,\ox,\bar\lm)$.
\end{definition}

Our subsequent goal is to establish relationships between the KKT strong regularity and full stability of local minimizers in composite optimization. We show below that these notions are actually equivalent under nondegeneracy. The result obtained below continues the line of equivalencies developed recently for various problems of constrained optimization in \cite{mnr,mrs,ms14} while being new for the composite optimization problems studied in the paper. To proceed, we consider the following {\em canonically perturbed} version $\Tilde{{\cal P}}_{\ow}(v_1,v_2)$ of problem (\ref{fcp1}) with parametric pairs $(v_1,v_2)\in\R^n\times\R^m$:
\begin{equation}\label{ccp}
\mbox{minimize }\;\ph_0(x,\ow)+\th(\Phi(x,\ow)+v_2)-\la v_1,x\ra\;\mbox{ subject to }\;x\in\R^n.
\end{equation}

The next lemma important in what follows reduces the study of full stability in the original optimization problem (\ref{fcp2}) to that in the canonically perturbed one  (\ref{ccp}) under nondegeneracy. Its proof is based on the presented criterion of full stability  in Theorem~\ref{ssooc} and allows us to deal with generalized equations of type (\ref{kkt11}) whose set-valued parts depend on parameters.

\begin{lemma}{\bf (full stability with respect to canonical perturbations).}\label{rsr2} Let $\ox$ be a feasible
solution to the composite optimization problem ${\cal P}(\ow,\ov)$ in {\rm(\ref{fcp2})} with some $\ow\in\R^d$ and $\ov$ from $(\rm{\ref{ov}})$ under the nondegeneracy condition ND from {\rm(\ref{fnond})}. Then $\ox$ is a fully stable local minimizer of ${\cal P}(\ow,\ov)$ if and only if it is a fully stable local minimizer of $\Tilde{{\cal P}}_{\ow}(\ov,0)$ in {\rm(\ref{ccp})}.
\end{lemma}
{\it Proof} It is easy to observe from the equivalent dual representation (\ref{2qcf2}) of the nondegeneracy condition (\ref{fnond}), obtained in \cite[Theorem~3.2]{ms16}, that ND for the canonically perturbed problem (\ref{ccp}) agrees with the one for the  fully perturbed problem (\ref{fcp2}). Suppose now that $\ox$ is a fully stable local minimizer of $\Tilde{{\cal P}}_{\ow}(\ov,0)$ and then apply Theorem~\ref{ssooc} to conclude that it is equivalent to the validity of the following inequality:
\begin{equation}\label{eq300}
\la u,\nabla^2_{xx}L_{\ow}(\ox,0,\bar\lm)u\ra>0\;\mbox{ for all }\;0\ne u\in {\cal S},
\end{equation}
where the subspace ${\cal S}$ is defined in (\ref{sssoc2}) and $L_{\ow}$ is the Lagrangian associated with problem (\ref{ccp}) given by $L_{\ow}(x,v_2,\lm)\colon =\ph_0(x,\ow)+\la \lm,\Phi(x,\ow)+v_2\ra$. Therefore we have $\nabla^2_{xx}L_{\ow}(\ox,0,\bar\lm)=\nabla^2_{xx}L(\ox,\ow,\bar\lm)$ with $L$ coming from (\ref{kkt}), which indeed tells us that $\ox$ is a fully stable local minimizer of ${\cal P}(\ow,\ov)$. The converse implication of the lemma is verified similarly. $\b$\vspace*{0.1in}

Now we are ready to establish the aforementioned relationships between full stability of local minimizers in composite optimization and strong regularity of the associated KKT systems.

\begin{theorem}{\bf (relationships between full stability and strong regularity in composite optimization).}\label{rsr-1} Let $\ox$ be a feasible solution to ${\cal P}(\ow,\ov)$ in {\rm(\ref{fcp2})} with some $\ow\in\R^d$ and $\ov=0$ from $(\rm{\ref{ov}})$. Assume that the qualification condition {\rm(\ref{2.9})} holds. Then for some $\gg>0$ the following are equivalent :

{\bf (i)} $\ox$ is a fully stable locally optimal solution to ${\cal P}(\ow,\ov)$ satisfying condition ND from {\rm(\ref{fnond})}.

{\bf (ii)} $\ox\in M_\gamma(\ow,\ov)$ and $(\ox,\bar\lm)$ is a strongly regular solution to {\rm(\ref{kkt11})}.
\end{theorem}
{\it Proof} We first verify implication (ii)$\Longrightarrow$(i). It has been well recognized (see, e.g., \cite[Theorem~2B.10]{dr} that strong regularity of the KKT system (\ref{kkt11}) at $(\ox,\bar\lm)$ is equivalent to the fact that the KKT system associated with the canonically perturbed problem (\ref{ccp}) and given by
\begin{eqnarray*}\label{p-kkt}
\left[\begin{array}{c}
v_1\\
v_2\\
\end{array}
\right]
\in\left[\begin{array}{c}
\nabla_x L(x,\ow,\lm)\\
-\Phi(x,\ow)\\
\end{array}
\right]+\left[\begin{array}{c}
0\\
(\sub\th)^{-1}(\lm)\\
\end{array}
\right]
\end{eqnarray*}
admits a Lipschitz continuous and  single-valued graphical localization around $(0,0,\ox,\bar\lm)\in\R^n\times\R^m\times\R^n\times\R^m$. Thus it results from Theorem~\ref{rsr} that the composite SSOSC from (\ref{eq300}) and the nondegeneracy condition ND for $\Tilde{{\cal P}}_{\ow}(\ov,0)$ are satisfied. As mentioned in the proof of Lemma~\ref{rsr2}, the nondegeneracy conditions ND for both problems $\Tilde{{\cal P}}_{\ow}(\ov,0)$ and ${\cal P}(\ow,\ov)$ are the same, and therefore Theorem~\ref{ssooc} says that $\ox$ is a fully stable local minimizer for $\Tilde{{\cal P}}_{\ow}(\ov,0)$. Employing Lemma~\ref{rsr2} tells us that $\ox$ is a fully stable local minimizer for the original problem ${\cal P}(\ow,\ov)$ as well, which justifies that (ii)$\Longrightarrow$(i). By similar arguments we verify the converse implication and thus complete the proof. $\b$\vspace*{0.05in}

As a by-product of the obtained equivalence and the characterization of full stability of local minimizers in Theorem~\ref{ssooc}, we get the composite SSOSC characterization of strong regularity for the associated KKT system (\ref{kkt11}). The results of this type for various problems of constrained optimization with ${\cal C}^2$-smooth data can be found in \cite{dr96,mnr,mrs,kk,bs} via appropriate SSOSC and nondegeneracy conditions. Note that, in contrast to full stability, the corresponding nondegeneracy condition is {\em necessary} for strong regularity. Some second-order characterizations of full stability {\em without nondegeneracy} have been recently established in \cite{mn3} for NLPs.\vspace*{0.05in}

The last part of this section is devoted to studying relationships between strong regularity in the sense of Definition~\ref{str} and the notion of strong Lipschitzian stability, which is a Lipschitzian version of Kojima's strong stability \cite{ko}. The concept of strong Lipschitzian stability was considered before only for problems of constrained optimization with ${\cal C}^2$-smooth data; \textcolor{red}{\xcancel{see \cite{mnr,bs}  Here we extend its} {see \cite{mnr,bs}.  Here we extend it}} to the general framework of composite optimization problems and then show that it is indeed equivalent to strong regularity of the corresponding KKT system. Note that relationships between strong regularity and strong stability were first studied in \cite{jmrt} for classical NLPs and then further developed for more general constrained problems in \cite{mnr,kk,bs,ms14}.

To proceed in our composite optimization setting, suppose without loss of generality that $\ov=0$ and say that the pair $(\xi(x,u),\Upsilon(x,u))$ with $u\in\R^q$, $\xi\colon\R^n\times\R^q\to\R$, and $\Upsilon\colon\R^n\times\R^q\to\R^m$ is a {\em ${\cal C}^2$-smooth parametrization}
of $(\ph_0(x,\ow),\Phi(x,\ow))$ in ${\cal P}(\ow,0)$ at $\ou\in\R^q$ provided that $\ph_0(x,\ow)=\xi(x,\ou)$ and $\Phi(x,\ow)=\Upsilon(x,\ou)$ for all $x\in\R^n$, where both functions $\xi$ and $\Upsilon$ are ${\cal C}^2$-smooth. Consider now the family of the parametric optimization problems given by
\begin{eqnarray*}
\hat{{\cal P}}(u):\quad\mbox{minimize }\;\xi(x,u)+\th(\Upsilon(x,u))\;\mbox{ subject to }\;x\in\R^n.
\end{eqnarray*}

\begin{definition}{\bf(strong Lipschitzian stability for composite optimization problems).}\label{sls}
A stationary point $\ox$ of problem ${\cal P}(\ow,0)$ in {\rm(\ref{fcp2})} is called {\sc strongly Lipschitz stable} with respect to the ${\cal C}^2$-smooth parametrization $(\xi(x,u),\Upsilon(x,u))$ of $(\ph_0(x,\ow),\Phi(x,\ow))$ in ${\cal P}(\ow,0)$ at $\ou\in\R^q$ if there are neighborhoods $U$ of $\ou$ and $O$ of $\ox$ such that for any $u\in U$ each problem $\hat{{\cal P}}(u)$ has a unique stationary point $x(u)\in O$ and the mapping $u\longmapsto x(u)$ is locally Lipschitzian around $\ou$. If it holds for any ${\cal C}^2$-smooth parameterization of $(\ph_0(x,\ow),\Phi(x,\ow))$ in ${\cal P}(\ow,0)$ at $\ou\in\R^q$, then the stationary point $\ox$ is called strongly Lipschitz stable.
\end{definition}

\begin{theorem}{\bf (equivalence between strong regularity and strong Lipschitzian stability for composite optimization problems).}\label{sls2}
Let $\ox$ be a feasible solution to the unperturbed problem ${\cal P}(\ow,\ov)$ in {\rm(\ref{fcp2})} with some $\ow\in\R^d$ and $\ov=0$ from $(\rm{\ref{ov}})$.  Assume further that the qualification condition {\rm(\ref{2.9})}
holds. Then the following are equivalent for some $\gg>0$:

{\bf (i)} $\ox$ is a Lipschitz stable local optimal minimizer of ${\cal P}(\ow,0)$ satisfying condition ND from {\rm(\ref{fnond})}.

{\bf (ii)} $\ox\in M_\gamma(\ow,\ov)$ and $(\ox,\bar\lm)$ is a strongly regular solution to {\rm(\ref{kkt11})}.
\end{theorem}
{\it Proof} Suppose that (i) holds. Since $(\ph_0(x,w)-\la x,v\ra,\Phi(x,w))$ is a ${\cal C}^2$-smooth parametrization of $(\ph_0(x,\ow),\Phi(x,\ow))$ in problem ${\cal P}(\ow,0)$ in {\rm(\ref{fcp2})} at the point $\ou:=(\ow,0)\in\R^d\times\R^n$, we find some neighborhoods $U$ of $\ou$ and $O$ of $\ox$ such that for any $u=(w,v)\in U$ there exists a unique stationary point $x(u)$ of $\hat{{\cal P}}(u)$ for which the mapping $u\longmapsto x(u)$ is Lipschitz continuous around $(\ou,\ox)$. This shows that the set-valued mapping
\begin{eqnarray*}
S(u):=\Big\{x\in\R^n\,\colon\;v\in\nabla_x\ph_0(x,w)+\nabla_x\Phi(x,w)^*\sub{\th}(\Phi(x,w))\Big\}
\end{eqnarray*}
admits a Lipschitzian single-valued graphical localization around $(\ou,\ox)$. Employing \cite[Theorem~3.4]{mrs}, we see that $\ox$ is a fully stable
locally optimal solution to problem ${\cal P}(\ow,0)$, which in turn yields the validity of (ii) due to Theorem~\ref{rsr-1}.

To prove the converse implication ${\rm(ii)}\Longrightarrow{\rm(i)}$, let $(\ox,\bar\lm)$ be a strongly regular solution to the KKT system {\rm(\ref{kkt11})}. This tells us that $\ox$ is a fully stable local minimizer of ${\cal P}(\ow,0)$ due to Theorem~\ref{rsr-1} and that the nondegeneracy condition ND is satisfied. Pick now an arbitrary ${\cal C}^2$-smooth parametrization $(\xi(x,u),\Upsilon(x,u))$ of $(\ph_0(x,\ow),\Phi(x,\ow))$ in ${\cal P}(\ow,0)$ at $\ou\in\R^q$,
which gives us the equalities $\nabla_x\ph_0(\ox,\ow)=\nabla_x\xi(\ox,\ou)$ and $\nabla_x\Phi(\ox,\ow)=\nabla_x\Upsilon(\ox,\ou)$ together with  those for the corresponding second-order derivatives. Therefore the composite SSOSC from (\ref{sssoc1}) is satisfied for problem $\hat{{\cal P}}(\ou)$, which in turn implies that $\ox$ is a fully stable local minimizer of problem $\hat{{\cal P}}(\ou)$. Employing now \cite[Theorem~3.4]{mrs}, we deduce that the set-valued mapping
\begin{eqnarray*}
S(u,v):=\Big\{x\in\R^n\,\colon\;v\in\nabla_x\xi(x,u)+\nabla_x\Upsilon(x,u)^*\sub{\th}(\Upsilon(x,u))\Big\}
\end{eqnarray*}
admits a Lipschitzian single-valued graphical localization around $(\ou,0,\ox)$. Defining $x(u):=S(u,0)$, conclude that it is a stationary point for problem $\hat{{\cal P}}(u)$ and that the mapping $u\longmapsto x(u)$ is locally Lipschitzian around $(\ou,\ox)$. This verifies (i) and completes the proof of theorem.
$\b$\vspace*{-0.2in}

\section{Lipschitzian Stability of Parametric Variational Systems}\label{Lip-PVS}

This section concerns Lipschitzian stability of solution maps to {\em parameterized generalized equations} (in Robinson's terminology) written in the form
\begin{eqnarray}\label{ge}
0\in f(x,w)+F(x,w)\;\mbox{ with }\;x\in\R^n\;\mbox{ and }\;w\in\R^d,
\end{eqnarray}
where $x$ is the decision variable and  $w\in\R^d$ stands for parameters, and where $f\colon\R^n\times\R^d\to\R^q$ is a single-valued {\em base} and $F\colon\R^n\times\R^d\tto\R^q$ is a set-valued {\em field}. We consider here more specific forms of (\ref{ge}), where $F$ is described by a  subdifferential mapping generated by an extended-real-valued function (called {\em potential}) involving a CPWL one. Due to their subdifferential/normal cone structure, such systems reflect certain variational properties as, e.g., the KKT system (\ref{kkt11}) associated with the composite optimization problem (\ref{fcp2}) under perturbations. Thus we use in what follows the term {\em parametric variational systems} (PVS) for generalized equations of this type.

Our first object to consider here is the PVS {\em solution map} given by
\begin{equation}\label{geq}
S(w):=\Big\{x\in\R^m\,\colon\;0\in f(x,w)+\sub\th(x)\Big\},
\end{equation}
where $\th\colon\R^m\to\oR$ is a parameter-independent CPWL function. We see that system (\ref{geq}) is the solution map to (\ref{ge}) with $F(x,w)=\partial\th(x)$, $q=n$, and $m=n$. In the particular case of $\th(x)=\dd(x;\O)$ the generalized equation in (\ref{geq}) amounts to the classical (parameterized) variational inequality, while in the general case of a convex function $\th$ in (\ref{geq}) such a system is called sometimes ``variational inequality of the second kind."

Our main goal in what follows is to study Lipschitzian stability of PVS in the sense that the solution map satisfies the so-called {\em Lipschitz-like/Aubin property} around the given point $(\ow,\ox)\in\gph S$: there are neighborhoods $U$ of $\ow$ and $O$ of $\ox$ and a number $\ell>0$ such that
\begin{equation}\label{lipl}
S(w_1)\cap O\subset S(w_2)+\ell\|w_1-w_2\|\B\;\mbox{ for all }\;w_1,w_2\in U,
\end{equation}
where $\B$ stands for the closed unit ball in $\R^n$. The infimum of all such moduli $\ell$ appearing in (\ref{lipl}) is called the {\em exact Lipschitzian bound } of $S$ around $(\ow,\ox)$ and is denoted by $\lip S(\ow,\ox)$.

We know from \cite[Theorem~5.7]{m93} and \cite[Theorem~9.40]{rw} that the latter property can be  characterized by the following {\em coderivative/Mordukhovich criterion}, where $S$ should be closed-graph near $(\ow,\ox)$, which will be proved for the mappings considered below:
\begin{eqnarray}\label{cod2}
\begin{array}{ll}
D^*S(\ow,\ox)(0)=\{0\}\;\mbox{ with}\\
\lip S(\ow,\ox)=\sup\Big\{\|a\|\;\colon a\in D^*S(\ow,\ox)(b),\|b\|\le 1\Big\}.
\end{array}
\end{eqnarray}

Based on (\ref{cod2}) and coderivative calculus rules, characterizations and sufficient conditions for the Lipschitz-like property of solution maps
to generalized equations (\ref{ge}) and their specifications were obtained in \cite{m06,m94} and other publications in terms of the coderivative of set-valued mappings therein. Now we derive explicit characterizations of this property for $S$ from (\ref{geq}) (and then for some other types of PVS) via the {\em initial data} of the potential.

Given $\th\in CPWL$ described in (\ref{theta}) and (\ref{dom}), recall \cite[Definition~5.6]{ms15} that the {\em affine independence constraint qualification} (AICQ) is satisfied for the generating vectors $a_i$ and $d_t$  therein indexed by $(i,j)\in K(\ox)\times I(\ox)$ if the vectors $\{(a_i,1)\in\R^n\times\R|\;i\in K(\ox)\}\cup\{(d_t,0)\in\R^n\times\R|\;t\in I(\ox)\}$ are linearly independent. For $(\ox,\ov)\in\gph\th$ and $u\in\dom\sub^2\th(\ox,\ov)$, define
\begin{eqnarray*}
\begin{array}{ll}
I_{0,1}(u):=&\Big\{i\in K(\ox)\,\colon\;\la a_i-a_j,u\ra=0\;\;\mbox{as}\;\;j\in J_1\Big\},\\
I_{>,1}(u):=&\Big\{i\in K(\ox)\,\colon\;\la a_i-a_j,u\ra>0\;\;\mbox{as}\;\;j\in J_1\Big\},\\
I_{0,2}(u):=&\Big\{t\in I(\ox)\,\colon\;\la d_t,u\ra=0\Big\},\;I_{>,2}(u):=\Big\{t\in I(\ox)\,\colon\;\la d_t,u\ra>0\},
\end{array}
\end{eqnarray*}
where $J_1=J_+(\ox,\ov_1)$ and $J_2=J_+(\ox,\ov_2)$ with $\ov=\ov_1+\ov_2$.

Next we present an explicit characterization of the Lipschitz-like property of the solution map (\ref{geq}) in the general case of $\th\in CPWL$ in (\ref{geq}) as well as under the validity of AICQ.

\begin{theorem}{\bf(characterizations of Lipschitzian stability for PVS with CPWL potentials).}\label{llps} Let $(\ow,\ox)\in\gph S$ for the mapping $S$ from {\rm (\ref{geq})}, let $f$ be strictly differentiable at $(\ox,\ow)$, and let $\nabla_w f(\ox,\ow)$ be surjective. Denoting $\ov:=-f(\ox,\ow)\in\sub\th(\ox)$, we have the following statements:

{\bf(i)} The solution map $S$ is Lipschitz-like around $(\ow,\ox)$ if and only if
\begin{eqnarray*}
\Big[-\nabla_x f(\ox,\ow)^*u\in{\cal F}_{\tiny\{P_1,Q_1\},\{P_2,Q_2\}},\;-u\in{\cal G}_{\tiny\{P_1,Q_1\},\{P_2,Q_2\}}\Big]\Longrightarrow u=0
\end{eqnarray*}
whenever $(P_1,Q_1,P_2,Q_2)\in{\cal A}$ with the set ${\cal A}$ defined in {\rm (\ref{eq092})}. Furthermore, the exact Lipschitzian bound of $S$ around $(\ox,\ow)$ is calculated by
\begin{eqnarray*}\label{elb1}
\begin{array}{ll}
\lip S(\ow,\ox)=\sup\Big\{&\|\nabla_w f(\ox,\ow)^*u\|\;\colon-u\in{\cal G}_{\tiny\{P_1,Q_1\},\{P_2,Q_2\}},\;y\in{\cal F}_{\tiny\{P_1,Q_1\},\{P_2,Q_2\}},\\
&(P_1,Q_1,P_2,Q_2)\in{\cal A},\;\|y+\nabla_x f(\ox,\ow)^*u\|\le 1\Big\}.
\end{array}
\end{eqnarray*}

{\bf(ii)} Let AICQ hold in addition to {\rm(i)}. Then the solution map $S$ is Lipschitz-like around $(\ow,\ox)$ if and only if the only $u\in\R^n$ satisfying
\begin{eqnarray*}
\begin{array}{rll}
-u&\in&\dom\sub^2\th(\ox,\ov),\\
-\nabla_x f(\ox,\ow)^*u&\in &\span\{a_i-a_j\,\colon\;i,j\in I_{0,1}(u)\}\\&+&\cone\{a_i-a_j\,\colon\;i\in I_{>,1}(u),\;j\in I_{0,1}(u)\}\\
&+&\span\{d_t\,\colon\;t\in I_{0,2}(u)\}+\cone\{d_i\,\colon\;t\in I_{>,2}(u)\}
\end{array}
\end{eqnarray*}
is $u=0$, where the index sets $I:=I(\ox)$, $K\colon =K(\ox)$, $J_1:=J_+(\ox,\ov_1)$, and $J_2:=J_+(\ox,\ov_2)$ are defined by {\rm (\ref{active2})} and {\rm (\ref{eq05})}, respectively. Furthermore, the exact Lipschitzian bound of $S$ around $(\ox,\ow)$ is calculated by
\begin{eqnarray*}
\begin{array}{ll}
\lip S(\ow,\ox)=\sup\Big\{&\|\nabla_w f(\ox,\ow)^*u\|\,\colon\; -u\in{\cal G}_{\tiny \{J_1,J_1\},\{J_2,J_2\}},\\&\|w+\nabla_x f(\ox,\ow)^*u\|\le 1,\\
& w\in\span\{a_i-a_j\,\colon\;i,j\in I_{0,1}(u)\}\\&+\cone\{a_i-a_j\,\colon\;i\in I_{>,1}(u),\;j\in I_{0,1}(u)\}\\
&+\span\{d_t\,\colon\;t\in I_{0,2}(u)\}+\cone\{d_i\,\colon\;t\in I_{>,2}(u)\}\Big \}.
\end{array}
\end{eqnarray*}
\end{theorem}
{\it Proof} Let us first verify the following two facts.\vspace{0.05in}\\
{\bf{Fact~1:}} {\em The graph of the solution map $S$ from {\rm(\ref{geq})} is closed.}\\[1ex]
To verify this, take a sequence $(w_k,x_k)\in\gph S$ with $(w_k,x_k)\to(\hat w,\hat x)$ as $k\to\infty$. It yields $-f(x_k,w_k)\in \sub\th(x_k)$, which means that
$$
\la-f(x_k,w_k),x-x_k\ra\le\th(x)-\th(x_k)\;\mbox{ for all }\;x\in\dom\th.
$$
Passing there to the limit as $k\to\infty$ and taking into account that $\th$ is continuous relative to its domain by \cite[Proposition~10.21]{rw}, we conclude that
$$
\la-f(\hat x,\hat w),x-\ox\ra\le\th(x)-\th(\hat x)\;\mbox{ whenever }\;x\in\dom\th.
$$
Thus we arrive at $-f(\hat x,\hat w)\in\sub\th(\hat x)$, which yields $(\hat w,\hat x)\in\gph S$.\vspace{0.05in}\\
{\bf{Fact~2:}} {\em The coderivative of the solution map $S$ is represented by}
\begin{equation}\label{cod1}
\begin{array}{ll}
D^*S(\ow,\ox)(b)=\Big\{a\in\R^m\,\colon&\exists\;u\in\dom\sub^2\th(\ox,\ov)\;\;\mbox{and}\;\;a=\nabla_w f(\ox,\ow)^*u,\\
&-\nabla_x f(\ox,\ow)^*u-b\in\sub^2\th(\ox,\ov)(u)\Big\}.
\end{array}
\end{equation}
This follows from \cite[Theorem~4.44]{m06} due to the assumptions made.

Now we proceed with verifying both assertions (i) and (ii) simultaneously by using the coderivative criterion with the exact bound formula in (\ref{cod2}). Observe that $u=0$ amounts to $\nabla_w f(\ox,\ow)^*u=0$ due to the surjectivity of $\nabla_w f(\ox,\ow)$. The necessary and sufficient conditions in (i) and (ii) come out directly from formula (\ref{2nd-val}) and \cite[Theorem~4.10]{ms15}, respectively. The exact bound formulas (i) and (ii) follow from (\ref{cod2}) applied to the solution mapping (\ref{geq}) and its coderivative representation in (\ref{cod1}). This completes the proof.$\b$

Note that Theorem~\ref{llps} extends the results of \cite[Theorem~5.3]{hmn} from the case of $\th=\dd_Z$, the indicator function of a convex polyhedron, to the case of a general CPWL function $\th\colon\R^m\to\oR$.

Next we consider solution maps of PVS with {\em fully amenable} potentials
\begin{equation}\label{geq2}
S(w):=\Big\{x\in\R^n\,\colon \;0\in f(x,w)+\sub_x(\th\circ\Phi)(x,w)\Big\},
\end{equation}
where $f$ is the same as in (\ref{geq}) while $\th\in CPWL$ and $\Phi\colon\R^n\times\R^d\to\R^m$ is a ${\cal C}^2$-smooth mapping around the reference point $(\ox,\ow)$. This case is significantly more involved in comparison with (\ref{geq}), particularly due to the {\em parameter-dependent} field in (\ref{geq2}). The following theorem gives us sufficient conditions for the Lipschitz-like property of the solution map (\ref{geq2}).

\begin{theorem}\label{llps2}{\bf(Lipschitzian stability for PVS with parameter dependent fields).} Let $(\ow,\ox)\in\gph S$ for $S$ from {\rm(\ref{geq2})}, let $f$ be strictly differentiable at $(\ox,\ow)$, and let $\ov:=-f(\ox,\ow)\in\sub_x(\th\circ\Phi)(\ox,\ow)$. Assume that the nondegeneracy condition ND from {\rm(\ref{fnond})} holds and that
\begin{equation}\label{imp1}
\left[
\begin{array}{ll}
0\in \nabla f(\ox,\ow)^*u+\Big(\nabla^2_{xx}\la\oy,\Phi\ra(\ox,\ow)u,\nabla^2_{xw}\la\oy,\Phi\ra(\ox,\ow)u\Big)\\
~~+\Big(\nabla_x\Phi(\ox,\ow),\nabla_w\Phi(\ox,\ow)\Big)^*\partial^2\th(\oz,\oy)(\nabla_x\Phi(\ox,\ow)u)
\end{array}
\right]\Longrightarrow u=0,
\end{equation}
where $\oy\in\sub\th(\Phi(\ox,\ow))$ is such that $\nabla_x\Phi(\ox,\ow)^*\oy=\ov$. Remembering that the second-order subdifferential of $\th$ is calculated in {\rm(\ref{2nd-val})}, we claim that $S$ is Lipschitz-like around $(\ow,\ox)$ if $a=0$ is the only vector of $\R^d$ for which the following conditions hold with some $u\in\R^n$:
\begin{eqnarray*}
\begin{array}{ll}
\nabla_x\Phi(\ox,\ow)u\in\dom\sub^2\th(\ox,\ov),\\(-\nabla_x f(\ox,\ow)^*u,a-\nabla_w f(\ox,\ow)^*u)\in\Big(\nabla^2_{xx}\la\oy,\Phi\ra(\ox,\ow)u,\nabla^2_{xw}\la\oy,\Phi\ra(\ox,\ow)u\Big)\\
\hspace{3.8cm}+\Big(\nabla_x\Phi(\ox,\ow),\nabla_w\Phi(\ox,\ow)\Big)^*\partial^2\th(\oz,\oy)(\nabla_x\Phi(\ox,\ow)u),
\end{array}
\end{eqnarray*}
\end{theorem}
{\it Proof} Similarly to Theorem~\ref{llps}, we proceed as follows:\vspace{0.05in}\\
{\bf{Fact~1:}} {\it The graph of $S$ from {\rm(\ref{geq2})} is locally closed around $(\ox,\ow)$.}\\[1ex]
To justify this fact, recall that condition ND implies the validity of (\ref{2.9}). Consider now neighborhoods $O$ of $\ox$ and $W$ of $\ow$ so that (\ref{2.9}) is fulfilled for any $(x,w)\in O\times W$. Select $\epsilon>0$ small to ensure that $\B_{\epsilon}(\ox,\ow)\subset O\times W$ and then check that the set $\gph S\cap \B_{\epsilon}(\ox,\ow)$ is closed. To see it, pick a sequence $(w_k,x_k)\in\gph S\cap\B_{\epsilon}(\ox,\ow)$ with $(w_k,x_k)\to(\hat w,\hat x)$ as $k\to\infty$.  By $-f(x_k,w_k)\in\sub_x(\th\circ\Phi)(x_k,w_k)=\nabla_x\Phi(x_k,w_k)^*\sub\th(z_k)$ with $z_k=\Phi(x_k,w_k)$ there are $p_k\in\sub\th(z_k)$ such that $-f(x_k,w_k)=\nabla_x\Phi(x_k,w_k)^*p_k$. Since $\{p_k\}$ is bounded due to (\ref{2.9}), suppose without loss of generality that $p_k\to\hat p$ for some $\hat p\in \sub\th(\hat z)$. This tells us that $-f(\hat x,\hat w)=\nabla_x \Phi(\hat x,\hat w)^*\hat p\in\Phi(\hat x,\hat w)^*\sub\th(\hat z)$,  which hence  yields $(\hat w,\hat x)\in\gph S$ and thus justifies the claimed fact.\vspace{0.1in}\\
{\bf{Fact~2:}} {\it We have the following upper estimate for the coderivative of $S$}:
\begin{eqnarray*}
\begin{array}{ll}
D^*S(\ow,\ox)(b)\subset\Big\{a\in\R^d\colon\exists\;u\in\R^n\;\mbox{with}\;(-b-\nabla_x f(\ox,\ow)^*u,a-\nabla_w f(\ox,\ow)^*u)\\ \hspace{4cm}\in\Big(\nabla^2_{xx}\la\oy,\Phi\ra(\ox,\ow)u,\nabla^2_{xw}\la\oy,\Phi\ra(\ox,\ow)u\Big)\\
\hspace{4cm}+\Big(\nabla_x\Phi(\ox,\ow),\nabla_w\Phi(\ox,\ow)\Big)^*\partial^2\th(\oz,\oy)(\nabla_x\Phi(\ox,\ow)u)\Big\}.
\end{array}
\end{eqnarray*}
Indeed, this follows from the result of \cite[Corollary~4.47]{m06} and the second-order subdifferential chain rule obtained above in \cite[Corollary~3.1]{ms16}.

To justify the sufficient conditions for Lipschitzian stability claimed in the theorem, we just need to implement the coderivative criterion (\ref{cod2}). $\b$

The last PVS we consider here is the solution map given by
\begin{equation}\label{geq3}
S(v):=\Big\{x\in\R^n\,\colon \;v\in f(x,\ow)+\sub_x(\th\circ\Phi)(x,\ow)\Big\}
\end{equation}
with the fixed basic parameter $\ow\in\R^d$. There are two crucial issues that distinguish (\ref{geq3}) from (\ref{geq2}): both base and field of (\ref{geq3}) are {\em parameter-independent} while the other (free) parameter $v\in\R^q$ enters the left-hand side of (\ref{geq3}). Such systems are known as {\em canonically perturbed} ones.

It is interesting to observe that (\ref{geq3}) appears from KKT conditions associated with {\em tilt perturbations} (cf.\ \cite{pr98} for the general format of unconstrained optimization with extended-real-valued objectives) of the original composite optimization problem (\ref{fcp1}), i.e., (\ref{fcp2}) with the fixed basic parameter $w=\ow$:
\begin{equation}\label{fcp22}
{\cal P}_{\ow}(v):\quad\mbox{minimize }\;\ph_0(x,\ow)+\theta(\Phi(x,\ow))-\la v,x\ra\;\mbox{ subject to }\;x\in\R^n.
\end{equation}
Indeed, the KKT system (\ref{kkt}) for (\ref{fcp22}) can be written in the form of (\ref{geq3}):
$$
v\in\nabla_x\ph_0(x,\ow)+\sub_x(\th\circ\Phi)(x,\ow).
$$

The next theorem gives us {\em necessary and sufficient} conditions for the Lipschitz-like property of (\ref{geq3}). Observe that in contrast to Theorem~\ref{llps2} we do not need to assume the validity of (\ref{imp1}). Furthermore, the surjectivity condition of Theorem~\ref{llps} is automatic in the setting of (\ref{geq3}).

\begin{theorem}{\bf(characterization of Lipschitzian stability of PVS with canonical perturbations).}\label{llps3} Let $(\ov,\ox)\in\gph S$ for the mapping $S$ from {\rm (\ref{geq3})}, where $f$ is strictly differentiable at $(\ox,\ow)$ with respect to $x$, $\th\in CPWL$, and condition ND from {\rm(\ref{fnond})} is satisfied. Denote $\op:=-f(\ox,\ow)\in\sub_x(\th\circ\Phi)(\ox,\ow)$ and take $\oy\in\sub\th(\Phi(\ox,\ow))$ such that $\nabla_x\Phi(\ox,\ow)^*\oy=\op$. Then the solution map $S$ is Lipschitz-like around $(\ov,\ox)$ if and only if $u=0$ is the only vector satisfying
\begin{equation}\label{char6}
0\in\nabla_x f(\ox,\ow)^*u+\nabla^2_{xx}\la\oy,\Phi\ra(\ox,\ow)u+\nabla_x\Phi(\ox,\ow)^*\sub^2\th(\ox,\op)(\nabla_x\Phi(\ox,\ow)u),
\end{equation}
where the second-order subdifferential of $\th$ is calculated in {\rm(\ref{2nd-val})}.
\end{theorem}
{\it Proof} Similar to Fact~1 in the proof of Theorem~\ref{llps2} we can clarify that the graph of $S$ is locally closed around $(\ov,\ox)$. To justify (\ref{char6}) as a characterization of the Lipschitz-like property of $S$ in {\rm(\ref{geq3})}, define $\Xi(x):=\sub_x(\th\circ\Phi)(x,\ow)$ and $g(v,x):=(x,v-f(x,\ow))$ and then get
$$
(v,x)\in\gph S\Longleftrightarrow g(v,x)\in\gph\Xi.
$$
It is easy to observe from the construction of $g$ that
$$\nabla g(\ov,\ox)=\left(\begin{array}{cc}
0&I_n \\
I_n&-\nabla_x f(\ox,\ow)^*\\
\end{array}\right),
$$
where $I_n$ stands for the $n\times n$ identity matrix. Thus the Jacobian $\nabla g(\ov,\ox)$ is of full rank. Appealing now to \cite[Theorem~1.17]{m06} tells us that $$N((\ov,\ox);\gph S)=\nabla g(\ov,\ox)^*N(g(\ov,\ox);\gph\Xi),$$ and thus we have the equivalence
$$
a\in D^* S(\ov,\ox)(b)\Longleftrightarrow c\in D^*\Xi(g(\ov,\ox))(d)\;\mbox{ with }\;\left\{\begin{array}{ll}
a=-d,\\
b=-(c+\nabla_x f(\ox,\ow)^*d).
\end{array}\right.
$$
Employing this together with \cite[Corollary~3.1]{ms16} gives us
$$
\begin{array}{ll}
D^*S(\ov,\ox)(b)=\Big\{a\;\colon &0\in b-\nabla_x f(\ox,\ow)^*a-\nabla^2_{xx}\la\oy,\Phi\ra(\ox,\ow)a \\&~~~+\nabla_x\Phi(\ox,\ow)^*\sub^2\th(\ox,\op)(-\nabla_x\Phi(\ox,\ow)a)\Big\}.
\end{array}
$$
To complete the proof, it remains to use the coderivative criterion (\ref{cod2}). $\b$\vspace*{0.05in}

The Lipschitz-like property of solution maps to canonically perturbed variational inequalities over polyhedral sets written as
\begin{equation}\label{geq5}
S(v)=\Big\{x\in\R^n\colon\;v\in f(x,\ow)+N(x;Z)\Big\}.
\end{equation}
has been addressed in \cite{dr96}. It is easy to see that (\ref{geq5}) corresponds to our setting in (\ref{geq3}) with the indicator function $\th=\dd_Z$ of the convex polyhedron $Z\subset\R^n$, $m=n$, and $\Phi(x,\ow)=x$. In this case our characterization (\ref{char6}) reduces to
\begin{equation}\label{fc1}
\Big [0\in\nabla_x f(\ox,\ow)^*u+\sub^2\th(\ox,\op)(u)\Big]\Longrightarrow u=0
\end{equation}
and by (\ref{2nd-val}) can be equivalently written entirely via the initial data
\begin{eqnarray*}
\Big[\nabla_x f(\ox,\ow)^*u\in{\cal F}_{\tiny\{P_1,Q_1\},\{P_2,Q_2\}}\;\mbox{ and }\;u\in{\cal G}_{\tiny\{P_1,Q_1\},\{P_2,Q_2\}}\Big]\Longrightarrow u=0
\end{eqnarray*}
for all $(P_1,Q_1,P_2,Q_2)\in{\cal A}$. Note that this characterization is much more efficient that the {\em critical face condition} obtained in \cite{dr96}, which involves closed faces of some polyhedral critical cone built upon the tangent cone to the convex polyhedron and is hard to be implemented. Certain specifications as well as extensions of the latter condition are derived in \cite{yy} for canonically perturbed affine variational inequalities, but some critical face expressions still remain therein. On the other hand, the results of \cite[Theorem~3]{dr96} and \cite[Theorem~3.9]{yy} establish the equivalence of the Lipschitz-like property of (\ref{geq5}) to Robinson's strong regularity, which postulates locally {\em single-valued} Lipschitzian behavior of the solution map. It is a challenging {\em open question} about the possibility to obtain such a result in the general CPWL framework of (\ref{geq3}).

To conclude this paper, recall the recent developments of \cite{mn15} on {\em full stability} of general parametric variational systems in the form
$$
v\in f(x,w)+\partial_x\ph(x,p),
$$
with parameters $(v,w)$, where $\ph$ belongs to a broad class of ``parametrically continuous prox-regular" functions \cite{lpr}. Such systems largely extend the PVS considered above, and the full stability property for them in the sense of \cite{mn15} generally {\em yields} the single-valuedness and local Lipschitz continuity of their solution maps. In the particular case of the variational inequalities (\ref{geq5}) the characterization of full stability from \cite[Theorem~4.8]{mn15} reads as
\begin{equation}\label{mn}
\la\nabla_x f(\ox,\ow)^*u,u\ra+\la q,u\ra>0\;\mbox{ for all }\;q\in\sub^2\th(\ox,\op)(u),\;u\ne 0,
\end{equation}
which implies (\ref{fc1}). However, the reverse implications fails in general as, e.g., for the case of $\th=0$ and $f\colon\R^2\to\R^2$ given by $f(x):=(x_1,-x_2)$ at $\ox=(0,0)$.\vspace*{-0.2in}

\section{Conclusions}
This paper continues the path started in \cite{ms16} about  various applications of the second-order subdifferential theory for CPWL functions recently developed in \cite{ms15}. The main result obtained in this vein reveals the equivalence between full stability and strong regularity for composite optimization problems under the validity of the nondegeneracy condition (\ref{fnond}). Recently, we obtained a counterpart of the reduction lemma \cite[Lemma~2E.4]{dr} for CPWL functions via a different approach, rooted in the established results in \cite{ms15}. This seems to be beneficial in the study of Lipschitzian stability for PVS (\ref{geq}). Employing this result, we plan to proceed with the study of the relationship between the Lipschitz-like property of (\ref{geq3}) and Robinson's strong regularity in the CPWL framework.\vspace*{-0.1in}

\begin{acknowledgements}
This research was partly supported by the National Science Foundation under grants DMS-1007132 and DMS-1512846 and by the Air Force Office of Scientific Research grant \#15RT0462. The authors are grateful to both anonymous referees for their valuable remarks, which helped us to improve the original presentation. \vspace*{-0.1in}
\end{acknowledgements}

\end{document}